\documentclass[10pt,oneside,reqno]{amsart}
\usepackage{srcltx}
\usepackage{pst-all}  
\usepackage{amssymb}
\usepackage{amsthm} 
\usepackage{amsmath,amsthm,amsopn,amstext,amsbsy}      
\usepackage{latexsym}
\usepackage{newlfont}
\usepackage{amsfonts}
\usepackage{mathrsfs}
\usepackage{comment}
\usepackage[affil-it]{authblk}
\newcommand{\R}{\mathbb{R}}

\newcommand{\Pb}{\mathbb{P}}
\newcommand{\Z}{\mathbb{Z}}
\newcommand{\N}{\mathbb{N}}
\newcommand{\be}{\begin{equation}}
\newcommand{\ee}{\end{equation}}
\def\non{\noindent}
\newtheorem{teo}{Theorem}

\newtheorem{prop}[teo]{Proposition}
\newtheorem{lem}[teo]{Lemma}

\newtheorem{obs}[teo]{Observation}
\title{Sub-Gaussian bound for the one-dimensional Bouchaud trap model}
\author{Manuel Cabezas}
\affil{Instituto Nacional de Matem\'atica Pura e Aplicada\\ Rio de Janeiro, Brazil}
\date{\today}
\email{mncabeza@mat.puc.cl}
\address{Instituto Nacional de Matem\'atica Pura e Aplicada\\Estrada Dona Castorina 110\\
  Rio de Janeiro, Brazil}
\begin{document}
\begin{abstract}
We establish a sub-Gaussian lower bound for the transition kernel of the one-dimensional, symmetric Bouchaud trap model, which provides a positive answer to the behavior predicted by Bertin and Bouchaud in \cite{bertin}. The proof rests on the Ray-Knight description of the local time of a one-dimensional Brownian motion. Using the same ideas we also prove the corresponding result for the FIN singular diffusion.
\end{abstract}
\maketitle
 \setlength{\baselineskip}{6mm}
\noindent {\footnotesize
}\\
\section{Introduction}\label{introduction}
The \textit{Bouchaud trap model} (BTM) is a random walk in a random medium. The medium represents a landcape composed of \textit{traps} which retain the walk for some amount of time.
In this paper we will consider the BTM taking values on the line $\Z$. Each $z\in\Z$ represents a trap of some depth $\tau_z>0$. For every fixed realization of the traps $\boldsymbol\tau:=(\tau_z)_{z\in\Z}$, $(X_t)_{t\geq0}$ is a continuous-time Markov chain with $X_0=0$ and whose jump rates are given by
 $$
c (x,y)
:=\left\{
  \begin{array}{ll}
    (2\tau_x)^{-1} \textrm{ if } |x-y|=1, \\
    0                 \textrm{                  otherwise. }
  \end{array}
\right.
 $$
 That is, $(X_t)_{t\geq0}$ is a random walk on $\Z$ which, each time it visits a site $z\in\Z$, it waits there an exponentially distributed time with mean $\tau_z$ and then jumps to $z-1,z+1$ with probability $1/2$ each. Hence $\tau_z$ should be regarded as the depth of the trap at $z$. The environment $\boldsymbol \tau$ is chosen at random, more precisely, $(\tau_z)_{z\in\Z}$ is an i.i.d.~ sequence of positive random variables.
The random walk $(X_t)_{t\geq0}$ defined in this way is the one-dimensional Bouchaud Trap Model.

When the environment is highly inhomogeneous (i.e, when the distribution of $\tau_0$ has heavy tails), the trapping mechanism becomes relevant even in the large scale behavior of the system, and the BTM displays some striking features, such as \textit{localization} (see \cite{fin99}), \textit{subdiffusivity} (see \cite{fin}) and \textit{aging} (see \cite{fin},\cite{cerny},\cite{bc}).
More precisely, \textit{localization} means that, for $\Pb$-a.e. realization of the environment $\boldsymbol\tau$, we have that
$$
\limsup_{t\to\infty}\max_{z\in\Z}\Pb(X_t=z\vert\boldsymbol\tau)>0,
$$
i.e, for arbitrarily large times, there is a site $z$ which carries a positive proportion of the distribution of $X_t$. Localization reflects the fact that, at all time scales, the depth of the deepest trap found by the BTM is of the same order of magnitude that the sum of the depths of all the traps visited by the BTM up to that time. 

Hence, the BTM is characterized by the effect of the deepest traps visited, which dominate the dynamics at all time scales. In agreement with this picture, the BTM has an anomalously  slow evolution. In fact, the \textit{subdiffusivity} of the BTM means that the typical displacement of $X$ at time $t$ is $o(t^{1/2})$.

Other interesting properties of the BTM can be found by analyzing some \textit{two-time} functions as
$f(t_w,t):=\Pb(X_{t_w+t}=X_{t_w})$. We should think that we are letting that the system \textit{ages} for a time $t_w$, and then we measure the probability that the system is at the ``initial position'' $X_{t_w}$ after a time $t$. Normally, one would expect that $f(t_w,t)$ depends only on $t$ (i.e, that the aging has no effect on the measurement). Nevertheless, for the BTM on a strongly inhomogeneous environment, we have that $f(t_w,t)$ depends on the ratio $t/t_w$ rather than on $t$. More precisely, as shown in \cite{fin}, we have that $$\lim_{\stackrel{t_w\to\infty}{ t/t_w\to\theta}}f(t_w,t)=F(\theta),$$ where $F$ is a non-trivial function of $\theta$. This dependence on the ratio $t/t_w$ for some two-time functions is usually referred as \textit{aging}. The aging exhibited by $f(t,t_w)$ is yet another consequence of the ever slower dynamics of the BTM.

The standard assumption to ensure an inhomogeneous environment is that the law of $\tau_0$ is in the domain of attraction of an $\alpha$-stable law, for some $\alpha\in(0,1)$. For the sake of simplicity, in this article we will make the stronger assumption that there exists $\alpha\in(0,1)$ such that
$$
\lim_{u\rightarrow\infty}u^\alpha\mathbb{P}(\tau_x\geq u)=1.
$$

A basic question of the model is to describe its transition kernel. In \cite{bertin}, E.M. Bertin and J.-P. Bouchaud predicted that the anomalously slow dynamics should be reflected on a non-Gaussian diffusion front. Moreover, they conjectured that the transition kernel has a decay given by a stretched exponential. More precisely, they claimed that $\Pb(|X_t|\geq x)$ should behave as $C\exp\left(-c\left(\frac{x}{t^{\gamma}}\right)^{1+\alpha}\right)$ for some positive constants $C,c$, where $\gamma=\alpha/(1+\alpha)$. Their prediction was supported by numerical simulations and non-rigorous arguments.
In this article we confirm the predicted behavior by providing the sub-Gaussian lower bound for the transition kernel. The corresponding upper bound was previously obtained by J. {\v{C}}ern{\'y} in \cite{cerny}.

We proceed to state the main result obtained in this article:
\begin{teo}\label{main}
There exists positive constants $C_1,c_1,C_2,c_2$ and $\epsilon_1$ such that
$$
 C_1\exp\left(-c_1\left(\frac{x}{t^{\gamma}}\right)^{1+\alpha}\right)\leq\Pb(|X_t|\geq x)\leq C_2\exp\left(-c_2\left(\frac{x}{t^{\gamma}}\right)^{1+\alpha}\right)
$$
for all $t\geq0$ and $x\geq 0$ such that $\frac{x}{\epsilon_1}\leq t$.
\end{teo}
\non As we have previously stated, the upper bound in Theorem \ref{main} has been already obtained in \cite[Lemma 3.2]{cerny}.

 The proof that we will present for the lower bound relies heavily on the fact that $(X_t)_{t\geq0}$ has a clearly identified scaling limit. This scaling limit is the \textit{Fontes-Isopi-Newman singular diffusion} (FIN) $(Z_t)_{t\geq0}$. It was discovered by  Fontes, Isopi and Newman in \cite{fin} and it should be thought of as a Brownian motion moving among a random environment of traps. More accurately, the FIN singular diffusion $Z$ is obtained as a speed-measure change of a
Brownian motion through a random, purely atomic measure $\rho$, where $\rho$ is the Stieltjes measure associated to an $\alpha$-stable subordinator.

It is a known fact that $\alpha$-stable subordinators are pure-jump processes. Hence, we can write
\begin{equation}\label{d:rho}
 \rho=\sum_{i\in\N} v_i \delta_{x_i}.
\end{equation}
Each atom $v_i\delta_{x_i}$ of $\rho$ plays the role of a trap located at $x_i$ and whose \textit{depth} is $v_i$. It can be shown that $\{x_i:i\in\N \}$ is dense in $\R$. Hence, the FIN diffusion is a Brownian motion moving among a random, dense set of traps $v_i\delta_{x_i},i\in\N$. 

As we have said, the FIN singular diffusion is the scaling limit of the one-dimensional Bouchaud trap model, more precisely, as proved in \cite[Theorem 4.1]{fin}, we have that
\begin{align*}
(\epsilon^{-1} X_{\epsilon^{1/\gamma}t})_{t\geq0} \stackrel{d}{\to} (Z_t)_{t\geq0} \quad \text{as } \epsilon \to 0, 
\end{align*}
 where $\stackrel{d}{\to}$ denotes convergence in distribution. Note that the scaling in the display above is sub-diffusive, in accordance with the subdiffusivity of the BTM.

The techniques used to prove Theorem \ref{main} can also be applied to establish the sub-Gaussian behavior of the transition kernel for the FIN singular diffusion, i.e, we will prove:
\begin{teo}\label{main2}
There exists positive constants $C_3,c_3,C_4$ and $c_4$ such that
$$
 C_3\exp\left(-c_3\left(\frac{x}{t^\gamma}\right)^{1+\alpha}\right)\leq\Pb(|Z_t|\geq x)\leq C_4\exp\left(-c_4\left(\frac{x}{t^\gamma}\right)^{1+\alpha}\right)
$$
for all $t\geq0$ and $x\geq 0$.
\end{teo}
\non Again, the upper bound of Theorem \ref{main2} was obtained by J. {\v{C}}ern{\'y} in \cite[Corollary 3.3]{cerny}.

At this point we would like to give some references concerning the BTM. The model was introduced by Bouchaud in \cite{weakergodicitybreaking} as a toy model to study metastability on some complex system such as spin glasses. In this case the state space of the walk was a large complete graph. Each vertex $x$ of the complete graph represented a metastable state of the system, and the depths where chosen as $\tau_x:=\exp(-\beta E_x)$, where $E_x$ represented the energy barrier that the system must overcome to leave that metastable state and $\beta$ is the inverse temperature. When the state space is $\mathbb{Z}^d, d \geq 2$, the BTM has a behavior completely different from the one-dimensional case, as shown by Ben Arous and \v{C}ern\'{y} in \cite{zetade}, and by Ben Arous, \v{C}ern\'{y} and Mountford in \cite{two} (see also \cite{Mourrat11}). In these papers it is shown that the scaling limit of the BTM on $\Z^d$ ($d\geq2$) is the \textit{fractional kinetics process} (FK),
which is a time-change of a $d$-dimensional Brownian motion through the inverse of an $\alpha$-stable subordinator.
In \cite{bbg} and \cite{bbg2} Ben Arous, Bovier and Gayrard obtained aging properties of the model on the complete graph. For an spectral characterization of aging see \cite{ BovierFaggionato05}. 
More on the one-dimensional BTM can be found in \cite{bc}, where it was shown that the model displays several aging regimes. The hydrodynamic behavior of one-dimensional BTM can be found in \cite{ JaraLandimTexeira11}. The response of the one-dimensional BTM to a drift is the subject of \cite{zindy}. In \cite{GantertMortersWachtel} a drift which decays to $0$ is introduced, and a phase transition in terms of the speed of decay of the drift is identified (see also \cite{Cabezas10}). 
 A study of the BTM on a wider class of graphs can be found
on \cite{universal}.
For a general account on the mathematical study of the Bouchaud trap model and the FIN singular diffusion, we refer to \cite{bcnotes}.

The organization of the paper is as follows: In Section \ref{s:strategyoftheproof} we give some preliminaries and we briefly sketch the main ideas of the proofs. The proofs of the main results are the subject of Section \ref{s:proofs}. Subsection \ref{ss:proofofmain2} contains the proof of Theorem \ref{main2} and Subsection \ref{ss:proofofmain1} deals with the proof of Theorem \ref{main}.

 In this article $C$ and $c$ will represent positive constants which might change their values from line to line.

\textbf{Acknowledgements}:
 The author would like to thank Alejandro Ramirez for fruitful discussions. This research was supported by a fellowship of the National Commission on Science and Technology of Chile (Conicyt)\#29100243.

\section{Preliminaries and strategy of the proof}\label{s:strategyoftheproof} 
\subsection{Preliminaries}\label{ss:preliminaries}
Here we will give the precise definition of the FIN singular diffusion and we will state the Ray-Knight Theorem.
We start with the definition of the FIN diffusion.

To define the FIN singular diffusion, first we recall the definition of a speed-measure changed Brownian motion.
Let $(B_t)_{t\geq0}$ be a standard one-dimensional Brownian motion starting at zero. Let $l(t,x)$ be a jointly-continuous version of its local time.
Given any locally finite measure $\mu$ on $\mathbb{R}$, denote
$$
\phi_{\mu}(s):=\int_{\mathbb{R}}l(s,y)\mu(dy),
$$
and its right-continuous generalized inverse by
$$
\psi_{\mu}(t):=\inf\{s>0:\phi_{\mu}(s)>t \}.
$$
The speed-measure change of $B$ with speed-measure $\mu$, $(B[\mu]_t)_{t\geq0}$, is defined as 
\be\label{speedmeasurechangedbrownianmotion}
B[\mu]_t:=B_{\psi_{\mu}(t)}.
\ee
Now, we proceed to define the random measure appearing on the definition of the FIN singular diffusion.
Let $(V_x)_{x\in\R}$ be a two sided, $\alpha$-stable subordinator independent of $B$.
 Let $\rho$ be the Lebesgue-Stieltjes measure associated with $V$, that is $\rho(a,b]:=V_b-V_a$.
The process $(Z_t)_{t\geq0}$ defined as $Z_s:=B[\rho]_s$ is the FIN singular diffusion.

\begin{obs}\label{scaleinvariance} The measure $\rho$ has scale invariance in the sense that $\lambda^{-1/\alpha}\rho(0,\lambda)$ is distributed as $\rho(0,1)$ for all $\lambda>0$. Also, the Brownian motion $B$ is scale invariant in the sense that $(\lambda^{-1/2}B_{\lambda t})_{t\geq 0}$ is distributed as $(B_t)_{t\geq 0}$. Those two facts imply that $Z$ is scale invariant in the sense that $(\lambda^{-\gamma}Z_{\lambda t})_{t\geq0}$ is distributed as $(Z_t)_{t\geq0}$ for all $\lambda>0$.
This fact reflects that the FIN singular diffusion has a subdiffusive behavior.
\end{obs}

Next we will recall the main toll of our proof, that is, the Ray Knight description of the local time of the Brownian motion.
Recall that $B$ is a standard one-dimensional Brownian motion started at the origin and $l(t,x)$ is its local time. For any $b\in\R$ let $\tau_{b}:=\inf\{t\geq0:B_t=b\}$. The Ray-Knight Theorem (\cite{ray} and \cite{knight}) states that:
\begin{prop}(Ray-Knight)
For each $a>0$, the stochastic process $(l(\tau_{-a},x): x\geq -a)$ is Markovian.
 Moreover $(l(\tau_{-a},x):-a\leq x \leq 0)$
is distributed as a squared Bessel process of dimension $d=2$ started at $0$ ($l(\tau_{-a},-a)=0$). Further, $(l(\tau_{-a},x):x\geq0)$ is distributed as a squared Bessel process of dimension $d=0$ started at $l(\tau_{-a},0)$ and killed at $0$.
\end{prop}

\subsection{Strategy of the proof}\label{ss:strategy}
In what follows, we will sketch the proof of Theorem \ref{main2}. The proof of Theorem \ref{main} follows the same line of reasoning but the technical details are slightly more complicated.

The main step in the proof of Theorem \ref{main2} is to find a lower bound for $\Pb(\min\{Z_s:s\in[0,t]\}\leq-x)$. That is equivalent to find a lower bound for $\Pb(H_{-x}\leq t)$, where, for any $b\in\R$, $H_b:=\inf\{t\geq0: Z_t=b\}$. The idea of the proof is to divide the interval $[-x,0)$ into $n$ sub-intervals $\left[-x\frac{i}{n},-x\frac{(i-1)}{n}\right),i=1,\dots,n$. Then we consider the events:
\begin{equation}\label{eq:firstsplitting}
 O_i:=\left\{\textrm{The time spent by } Z \textrm{ in } \left[-x\frac{i}{n},-x\frac{(i-1)}{n}\right) \textrm{ up to } H_{-x} \textrm{ is less than } \frac{t}{n}\right\}, i=1,\dots,n
\end{equation}
 Clearly, we have that the intersection of those events implies that the time spent by $Z$ in the negative axis before hitting $-x$ is less than $t$. If we ignore for the moment the problem of controlling the time spent in the positive axis, we just have to control the probability of the intersection $\cap_{i=1}^n O_i$. At this point, we would like to use the product rule to compute the probability of such intersection. Unfortunately, we are dealing with events which are not independent, and we cannot use the product rule. Nevertheless, one can proceed by using a similar strategy, and obtain independence by virtue of the Ray-Knight Theorem.

First, note that, for $i=1,\dots,n$, the time spent by $Z$ in the interval $[-x\frac{i}{n},-x\frac{(i-1)}{n})$ up to $H_{-x}$ can be written as
\begin{equation}
\Delta_i:=\int_{[-x\frac{i}{n},-x\frac{(i-1)}{n})}l(\tau_{-x},u)\rho(du), 
\end{equation}
where $\rho$ is the measure appearing in the definition of the FIN singular diffusion and $l(\cdot,\cdot)$ is the Brownian local time.
Hence the events in \eqref{eq:firstsplitting} can be written as
$$
O_i=\left\{\Delta_i\leq \frac{t}{n} \right\}, i=1,\dots,n.
$$

  The first step will be to deal with the dependence between $O_1$ and $(O_i)_{i=2,\dots,n}$.
Since $\rho$ is defined as the Lebesgue-Stieltjes measure associated to a L\'{e}vy process (more precisely, an $\alpha$-stable subordinator), we have that $\rho$ has independent increments, in the sense that, for any sequence $(J_i)_{i\in \N}$ of disjoint Borel subsets of $\R$, we have that the family of random variables $(\rho(J_i))_{i\in\N}$ is independent. Moreover, by the Ray-Knight Theorem, we know that $(l(\tau_{-x},u))_{u\geq0}$ is a Markov process. The key step in our reasoning is that, by virtue of the independence of increments of $\rho$ and the Markovianity of $l(\tau_{-x},u)$, we have that $\Delta_1$ depends of $\Delta_i, i=2,\dots,n$ only through the value of $l(\tau_{-x},-x/n)$. Hence, $O_1$ depends of $O_i,i=1,\dots,n$ only through $l(\tau_{-x},-x/n)$. As we will see in Section \ref{s:proofs}, it is not hard to deal with that kind dependence. 

The same argument can be applied to deal with the dependence between $O_j$ and $O_i,i=j+1,\dots,n$, for all $j=1,\dots,n-1$. In this way we will obtain a formula for the probability of the intersection $\cap_{i=1}^nO_i$ which will play the role of the product rule. 
  
Finally, thanks to the scale invariance of the FIN singular diffusion one can choose the number of sub-intervals $n$ in such a way that the probability of the events $O_i$ is a constant $c$ which does not depend on $x,t$ or $n$. The proper choice of $n$ will be $n=(\frac{x}{t^\gamma})^{1+\alpha}$. Hence, by using our product rule-type formula, we will get that the probability of $\left\{\min\{Z_s:s\in[0,t]\}\leq-x\right\}$ will be bounded below by $c^{(\frac{x}{t^\gamma})^{1+\alpha}}$ and that is the type of result that we want.

\begin{section}{Proofs of Theorems \ref{main} and \ref{main2}}\label{s:proofs}

\begin{subsection}{Proof of Theorem \ref{main2}}\label{ss:proofofmain2}
We will prove only the lower bound, since the upper bound is proved in \cite[Lemma 3.2]{cerny}.
 \non Thanks to the scale invariance of $Z$ (see observation \ref{scaleinvariance}), to prove the lower bound in Theorem \ref{main2} it is enough to show that there are positive constants $C,c>0$ such that
 $$
 \Pb(|Z_1|\geq x)\geq C\exp(-c x^{1+\alpha}) \textrm{ for all } x\geq 0.
 $$
\noindent Moreover, both $\Pb(|Z_1|\geq x)$ and $C\exp(-c x^{1+\alpha})$ are decreasing in $x$. Hence it will suffice to show that there are positive constants $C,c>0$ such that
$$
\Pb(|Z_1|\geq n^{1/(1+\alpha)})\geq C\exp(-c n)\textrm{ for all }n\in\N.
$$
The main step in the proof of Theorem \ref{main2} is the following lemma:
\begin{lem}\label{hittingtimesfin}
There are positive constants $C,c>0$ such that
$$
\Pb\left(\min\{Z_s:s\in[0,1]\}\leq -n^{1/(1+\alpha)}\right)\geq C\exp(-c n)\textrm{ for all }n\in\N.
$$
\end{lem}
\begin{proof}
\noindent
Note that the event $\left\{\min\{Z_s:s\in[0,1]\}\leq -n^{1/(1+\alpha)}\right\}$ is equivalent to $\left\{H_{-n^{1/(1+\alpha)}}\leq1\right\}$. We recall that $\rho$ is the random measure appearing in display \eqref{d:rho}, $B$ is the Brownian motion used in the construction of $Z$ and $l(\cdot,\cdot)$ is its local time. We have that
 $$
  H_{-n^{1/(1+\alpha)}}=\int_{-n^{1/(1+\alpha)}}^{\infty}l(\tau_{-n^{1/(1+\alpha)}},u)\rho(du).
  $$
  We will divide $[-n^{1/(1+\alpha)},0)$ into $n$ subintervals $I_i:=[-in^{-\gamma},-(i-1)n^{-\gamma}), i=1,\dots,n$ (recall that $\gamma=\alpha/(1+\alpha)$).
Let us define the events $A_i, i=1,\dots,n-1$ as
$$
A_i:=\left\{\int_{I_i}l(\tau_{-n^{1/(1+\alpha)}},u)\rho(du)\leq \frac{1}{n}\right\}.
$$
 Also define
$$
A_n:=\left\{\int_{I_n}l(\tau_{-n^{1/(1+\alpha)}},u)\rho(du)\leq \frac{1}{2n}\right\}
$$
and
$$
A_0:=\left\{\int_{[0,\infty)}l(\tau_{-n^{1/(1+\alpha)}},u)\rho(du)\leq \frac{1}{2n} \right\}.
$$
Clearly
\be\label{pre-spliting}
\bigcap_{i=0}^{n}A_i\subset\{H_{-n^{1/(1+\alpha)}}\leq 1\}.
\ee
At this point we would like to use the product rule on the display above to obtain a lower bound for $\Pb(H_{-n^{1/(1+\alpha)}}\leq1)$, but the events $A_i, i =0,\dots,n$ are not independent.
 Nevertheless, we will harness the Markovianity of $l(\tau_{-n^{1/(1+\alpha)}},u)$ (guaranteed by the Ray-Knight Theorem) and the independence of increments of $\rho$ to obtain a formula for $\Pb(\cap_{i=0}^nA_i)$ which will play the role of a product rule.
We will start by dealing with the dependence between $A_0$ and $\cap_{i=1}^nA_i$.

As we have explained in Section \ref{ss:strategy}, the measure $\rho$ has independent increments.
The key observation here is that, by virtue of the independence of increments of $\rho$ and by the Markovianity of $(l(\tau_{-n^{1/(1+\alpha)}},u))_{u\geq-n^{1/(1+\alpha)}}$ we have that the random variable $\int_{[0,\infty)}l(\tau_{-n^{1/(1+\alpha)}},u)\rho(du)$ depends of the random variables $\int_{I_i}l(\tau_{-n^{1/(1+\alpha)}},u)\rho(du), i=1,\dots,n$ only trough the value of $l(\tau_{-n^{1/(1+\alpha)}},0)$. In particular, we have that
$A_0$ depends of $A_i,i=1,\dots,n$ only through the value of $l(\tau_{-n^{1/(1+\alpha)}},0)$.

 To get rid of this dependence we will condition on the event that $l(\tau_{-n^{1/(1+\alpha)}},0)$ is below a certain fixed level $n^{-\gamma}$. That is, we introduce the event
$$
L_1:=\left\{l(\tau_{-n^{1/(1+\alpha)}},0)\leq n^{-\gamma}\right\}.
$$
The strategy will be that, when conditioned on the event $L_1$, the event $A_0$ is stochastically dominated by $\{\int_{[0,\infty)} Y_u\rho(du)\leq \frac{1}{2n}\}$, where $(Y_u)_{u\geq0}$ is a squared Bessel process of dimension $d=0$ started at $n^{-\gamma}$ ($Y_0=n^{-\gamma}$) and independent of $\rho$. In this way we will get rid of the dependence between $A_0$ and $l(\tau_{-n^{1/(1+\alpha)}},0)$.

We start by writing
\be\label{firstsplit}
\Pb\left(\cap_{i=0}^{n}A_i\right)\geq\Pb\left(A_0\cap L_1\cap_{i=1}^n A_i\right)
=\Pb(A_0\vert L_1\cap_{i=1}^n A_i)\cdot\Pb(L_1\cap_{i=1}^n A_i).
\ee
On the other hand, as we have said, the Markovianity of $l(\tau_{-n^{1/(1+\alpha)}},\cdot)$ and the independence of increments of $\rho$ imply that $A_0$ depends of $L_1\cap_{i=1}^nA_i$ only through the value of $l(\tau_{-n^{1/(1+\alpha)}},0)$. Hence we can write
\be\label{split}
\Pb(A_0\vert L_1\cap_{i=1}^n A_i)=\int_0^{n^{-\gamma}}\Pb(A_0\vert l(\tau_{-n^{1/(1+\alpha)}},0)=y)\nu^0(dy),
\ee
where $\nu^0$ is the distribution of $l(\tau_{-n^{1/(1+\alpha)}},0)$ conditioned on $L_1\cap_{i=1}^n A_i$.
Let us make a comment about the meaning of the conditional probabilities $\Pb(A_0\vert l(\tau_{-n^{1/(1+\alpha)}},0)=y)$ appearing in \eqref{split}.
One can use any regular version of the conditional probabilities of $(l(\tau_{-n^{1/(1+\alpha)}},u))_{u\geq0}$ given $l(\tau_{-n^{1/(1+\alpha)}},0)$ to define the function $y\mapsto\Pb(A_0\vert l(\tau_{-n^{1/(1+\alpha)}},0)=y)$ up to a set of $\nu^0$-measure zero. Nevertheless, as we will see, it will be convenient to define $\Pb(A_0\vert l(\tau_{-n^{1/(1+\alpha)}},0)=y)$ for all $y\geq0$. To do so, we use the fact that $(l(\tau_{-n^{1/(1+\alpha)}},u))_{u\geq0}$ is a squared Bessel process, and Bessel processes can be defined as starting from any $y\geq0$. Hence, for all $y\geq0$ we define $\Pb(A_0\vert l(\tau_{-n^{1/(1+\alpha)}},0)=y)$ using the law of a $0$-dimensional squared Bessel process started at $y$.

Now, we will show that, for all $y\in[0,n^{-\gamma}]$, we have that
\be\label{stochasticdomination}
\Pb(A_0\vert l(\tau_{-n^{1/(1+\alpha)}},0)=y)\geq\Pb(A_0\vert l(\tau_{-n^{1/(1+\alpha)}},0)=n^{-\gamma}).
\ee
We begin by showing that, for all $y\in[0,n^{-\gamma}]$, we can couple two $0$-dimensional squared Bessel processes, $(X^1_t)_{t\geq0}$, $(X^2_t)_{t\geq0}$ with $X^1_0=y$, $X^2_0=n^{-\gamma}$ and such that $\Pb(X^1_s\leq X^2_s\textrm{ for all }s\geq0)=1$. To see this, let $X^1$, $X^3$ be two independent $0$-dimensional squared Bessel processes with $X^1_0=y$ and $X^3_0=n^{-\gamma}$. Define
\be
X^2_t:=\left\{\begin{array}{ll}
X^3_t :t\leq\inf\{s\geq0:X^1_s=X^3_s\}\\
X^1_t :t>\inf\{s\geq0:X^1_s=X^3_s\}
\end{array}
\right.
\ee
That is, $X^2$ evolves independently of $X^1$ until the first time they meet, after which $X^2$ follows the same trajectory of $X^1$.
Then, the Markovianity of the processes involved implies that $X^2$ is distributed as a $0$-dimensional squared Bessel process, and by the continuity of paths of Bessel processes we have that $X^1_s\leq X^2_s$ for all $s\geq0$. In particular we have that for any measure $\mu$ on $\R$
$$
\int_0^\infty X^1_u\mu(du)\leq\int_0^\infty X^2_u\mu(du).
$$
Recalling the definition of the event $A_0$, we see that the previous display implies \eqref{stochasticdomination}.
On the other hand, displays \eqref{stochasticdomination} plus \eqref{split} imply that
$$
\Pb(A_0\vert L_1\cap_{i=1}^nA_i)\geq \Pb(A_0\vert l(\tau_{-n^{1/(1+\alpha)}},0)=n^{-\gamma}).
$$
Hence, using \eqref{firstsplit} we get that
\be\label{eq:A_0}
\Pb\left(\cap_{i=0}^nA_i\right)
\geq\Pb(L_1\cap_{i=1}^nA_i)\Pb(A_0\vert l(\tau_{-n^{1/(1+\alpha)}},0)=n^{-\gamma})
\ee
In this way we have been able to deal with the dependence between $A_0$ and $A_i,i=1,\dots,n$.

The next step is to get rid of the dependence between $A_1$ and $A_i,i=2,\dots,n$. We will proceed analogously. Let
$$L_2:=\left\{l(\tau_{-n^{1/(1+\alpha)}},-n^{-\gamma})\leq n^{-\gamma}\right\}.$$
We have that
$$
\Pb\left(L_1\cap_{i=1}^{n}A_i\right)\geq\Pb\left(A_1\cap L_1 \cap L_2\cap_{i=2}^n A_i\right)
=\Pb(A_1\cap L_1\vert L_2\cap_{i=2}^n A_i)\cdot\Pb(L_2\cap_{i=2}^n A_i).
$$
Again, by the Markovianity of $l(\tau_{-n^{1/(1+\alpha)}},\cdot)$ and the independence of increments of $\rho$, we have that $A_1\cap L_1$ depends on $L_2\cap_{i=2}^nA_i$ only through the value of $l(\tau_{-n^{1/(1+\alpha)}},-n^{-\gamma})$. That is
$$
\Pb(A_1\cap L_1\vert L_2\cap_{i=2}^n A_i)=\int_0^{n^{-\gamma}}\Pb(A_1\cap L_1\vert l(\tau_{-n^{1/(1+\alpha)}},-n^{-\gamma})=y)\nu^1(dy),
$$
where $\nu^1$ is the distribution of $l(\tau_{-n^{1/(1+\alpha)}},-n^{-\gamma})$ conditioned on $L_2\cap_{i=2}^n A_i$.
We can use the same type of coupling that we have used to obtain \eqref{stochasticdomination} to see that
$$
\Pb(A_1\cap L_1\vert l(\tau_{-n^{1/(1+\alpha)}},-n^{-\gamma})=y)\geq\Pb(A_1\cap L_1\vert l(\tau_{-n^{1/(1+\alpha)}},-n^{-\gamma})=n^{-\gamma}).
$$
for all $y\leq n^{-\gamma}$. Hence, the same reasoning used to obtain \eqref{eq:A_0} gives us
\begin{equation}\label{e:A_1}
\Pb(L_1\cap_{i=1}^nA_i)\geq\Pb(L_2\cap_{i=2}^nA_i)\Pb(A_1\cap L_1\vert l(\tau_{-n^{1/(1+\alpha)}},-n^{-\gamma})=n^{-\gamma}),
\end{equation}
which allows us to deal with the dependence between $A_1$ and $A_i,i=2,\dots,n$.

By analogous arguments, we get that, for all $j=1,\dots,n-1$
\begin{equation}\label{e:A_j}
\Pb(L_j \cap_{i=j}^nA_i)\geq\Pb(L_{j+1}\cap_{i=j+1}^nA_i)\Pb(A_j\cap L_j\vert l(\tau_{-n^{1/(1+\alpha)}},-jn^{-\gamma})=n^{-\gamma})
\end{equation}
where
$$
L_j:=\left\{l(\tau_{-n^{1/(1+\alpha)}},-(j-1)n^{-\gamma})\leq n^{-\gamma}\right\}.
$$
Putting together displays \eqref{eq:A_0}, \eqref{e:A_1} and \eqref{e:A_j}, we get that
\be\label{keyestimate}
\begin{split}
\Pb(\cap_{i=0}^nA_i)\\
  \geq\Pb\left(A_n\cap L_n\right)\Pb\left( A_0\;\middle\vert\;l(\tau_{-n^{1/(1+\alpha)}},0)= n^{-\gamma}\right)\prod_{i=1}^{n-1}&\Pb\left(A_i\cap L_i\;\middle\vert\;l(\tau_{-n^{1/(1+\alpha)}},-in^{-\gamma})= n^{-\gamma}\right).
 \end{split}
\ee
The formula above is the sort of product rule we were after.

Next, note that using the spatial homogeneity of the random measure $\rho$, display (\ref{keyestimate}) yields
\be\label{spliting}
\begin{split}
\Pb(\cap_{i=0}^nA_i)\\
\geq\Pb(A_n\cap L_n)\Pb\left( A_0\;\middle\vert\;l(\tau_{-n^{1/(1+\alpha)}},0)= n^{-\gamma}\right)&(\Pb\left(A_1\cap L_1\;\middle\vert\; l(\tau_{-n^{1/(1+\alpha)}},-n^{-\gamma})= n^{-\gamma}\right))^{n-1}.
\end{split}
\ee

To get Lemma \ref{hittingtimesfin}, it only remains to show that the quantities appearing in the display above (i.e, $\Pb(A_n\cap L_n), \Pb\left( A_0\;\middle\vert\;l(\tau_{-n^{1/(1+\alpha)}},0)= n^{-\gamma}\right)$ and $\Pb\left(A_1\cap L_1\;\middle\vert\; l(\tau_{-n^{1/(1+\alpha)}},-n^{-\gamma})= n^{-\gamma}\right)$) are positive and do not depend on $n$. We will first show that those quantities are independent of $n$, which will follow from the scale invariance of the measure $\rho$ and the scale invariance of the local time of the Brownian motion. We start with $\Pb\left(A_1\cap L_1\;\middle\vert\; l(\tau_{-n^{1/(1+\alpha)}},-n^{-\gamma})= n^{-\gamma}\right)$. We have that $\Pb\left(A_1\cap L_1\;\middle\vert\; l(\tau_{-n^{1/(1+\alpha)}},-n^{-\gamma})= n^{-\gamma}\right)$ equals
$$
\Pb\left(\int_{-n^{-\gamma}}^0l(\tau_{-n^{1/(1+\alpha)}},u)\rho(du)\leq \frac{1}{n}; l(\tau_{-n^{1/(1+\alpha)}},0)\leq n^{-\gamma}\;\middle\vert\; l (\tau_{-n^{1/(1+\alpha)}},-n^{-\gamma})=n^{-\gamma}\right).
$$
Recall that, if $(Y_t)_{t\geq 0}$ is a squared Bessel process started at some $a\geq0$, then, for each $\lambda\geq0$, $(\lambda Y_{\lambda^{-1} t})_{t\geq0}$ is a squared Bessel process started at $\lambda a$. Hence, by choosing $\lambda=n^{-\gamma}$, we have that the last display equals
$$
\Pb\left( \int_{-n^{-\gamma}}^0n^{-\gamma}l(\tau_{-n^{1/(1+\alpha)}},un^{\gamma})\rho(du)\leq \frac{1}{n} ; l(\tau_{-n^{1/(1+\alpha)}},0)\leq 1\;\middle\vert\; l(\tau_{-n^{1/(1+\alpha)}},-1)=1\right).
$$
Using the scale invariance of the measure $\rho$ we obtain that the previous display equals
$$
\Pb\left( \int_{-n^{-\gamma}}^0l(\tau_{-n^{1/(1+\alpha)}},un^{\gamma})\rho(n^{-\gamma}du)\leq 1 ; l(\tau_{-n^{1/(1+\alpha)}},0)\leq 1 \;\middle\vert\;l(\tau_{-n^{1/(1+\alpha)}},-1)=1\right).
$$
Finally, let us perform a change of variables inside the integral in the display above to get
\be\label{uno}
\Pb\left(\int_{-1}^0l(\tau_{-n^{1/(1+\alpha)}},u)\rho(du)\leq 1; l(\tau_{-n^{1/(1+\alpha)}},0)\leq 1 \;\middle\vert\; l(\tau_{-n^{1/(1+\alpha)}},-1)=1 \right).
\ee
Hence $\Pb\left(A_1\cap L_1 \; \middle \vert \; l(\tau_{-n^{1/(1+\alpha)}},-n^{-\gamma})=n^{-\gamma}\right)$ is equal to (\ref{uno}), which clearly does not depend on $n$.

Similar arguments can be used to show that
\be\label{dos}
\Pb(A_n\cap L_n)=\Pb\left(\int_{-1}^0l(\tau_{-1},u)\rho(du)\leq \frac{1}{2}; l(\tau_{-1},0)\leq 1\right)
\ee
and
\be\label{tres}
\Pb\left( A_0\;\middle\vert\;l(\tau_{-n^{1/(1+\alpha)}},0)=n^{-\gamma}\right)=\Pb\left(\int_{0}^\infty l(\tau_{-n^{1/(1+\alpha)}},u)\rho(du)\leq \frac{1}{2}\;\middle\vert\;l(\tau_{-n^{1/(1+\alpha)}},0)= 1\right),
\ee
which shows that the last two expressions are independent of $n$.

  It only remains to show that the quantities appearing on displays (\ref{uno}), (\ref{dos}) and (\ref{tres}) are non-zero.
Using the independence between the Brownian motion $B$ and the random measure $\rho$, we have that display \eqref{uno} is bigger or equal than
$$
\Pb\left(\sup_{u\in[-1,0]}l(\tau_{-n^{1/(1+\alpha)}},u)\leq \delta^{-1}; l(\tau_{-n^{1/(1+\alpha)}},0)\leq 1 \;\middle\vert\; l(\tau_{-n^{1/(1+\alpha)}},-1)=1 \right)\cdot \Pb\left(\rho(-1,0]\leq \delta\right)
$$
for any $\delta\geq0$.
We can see that the above display is positive for $\delta<1$ by noticing that the first factor in the product above is positive for $\delta<1$ (this follows from standard properties of Bessel processes) and, on the other hand, we have that $\Pb(\rho[-1,0)\leq\delta)>0$ for each $\delta>0$. An analogous argument gives that the quantity appearing in \eqref{dos} is positive. To see the positivity of the quantity appearing in \eqref{tres} we need also to use the fact that, for any $t>0$, a $0$-dimensional Bessel process hits the origin before time $t$ with positive probability.
\end{proof}
We have managed to show Lemma \ref{hittingtimesfin}, which provides the sub-Gaussian lower bound for $\Pb(\min\{Z_s:s\in[0,1]\}\geq-n^{1/(1+\alpha)})$. To prove Lemma \ref{lemmafin} we need to obtain the same type of bound for $Z_1$ instead of $\min\{Z_s:s\in[0,1]\}$. The standard argument of applying symmetry of the process after the stopping time $H_{-n^{1/(1+\alpha)}}$ to obtain $\Pb(Z_1\leq -n^{1/(1+\alpha)})=\frac{1}{2}\Pb(\min\{Z_s:s\in[0,1]\}\leq-n^{1/(1+\alpha)})$ cannot be applied in our case, because the annealed FIN singular diffusion is not Markovian (In particular, the event $\{\min\{Z_s:s\in[0,1]\}\leq -n^{1/(1+\alpha)}\}$ should be positively correlated with the absence of very deep traps on $[-n^{1/(1+\alpha)},0]$).
 To overcome this problem we will need the following lemma, in the argument we follow the reasoning used in \cite[Lemma 3.2]{cerny}.
\begin{lem}\label{fromminzetatozeta}
 There exists positive constants $C,c>0$ such that
 $$
\Pb\left(\int_{-n^{1/(1+\alpha)}}^0 l(\tau_{-n^{1/(1+\alpha)}},u)\rho(du)\geq 1; Z_s \leq n^{-\gamma} \textrm{ for all } s \leq H_{-n^{1/(1+\alpha)}} \right)\geq C\exp(-cn)
$$
for all $n\in\N$.

\end{lem}
\begin{proof}
\noindent
We divide the interval $[-n^{1/(1+\alpha)},0)$ into $n$ subintervals $I_i:=[-in^{-\gamma},-(i-1)n^{-\gamma})$, $i=1,\dots,n$ (recall that $\gamma=\alpha/(1+\alpha)$). For $i=1,\dots,n$ we define
$$
T_i:=\int_{I_i} (l(\tau_{-in^{-\gamma}},u)-l(\tau_{-(i-1)n^{-\gamma}},u))\rho(du).
$$
That is, $T_i$ is the time spent by $Z$ on the interval $I_i$ between times $H_{-(i-1)n^{-\gamma}}$ and $H_{-in^{-\gamma}}$. Also define, for $i=1,\dots,n$
$$
\theta_i:=\min\left\{s\geq H_{-(i-1)n^{-\gamma}}: Z_s\in\left\{-in^{-\gamma},-(i-2)n^{-\gamma}\right\} \right\}
$$
\noindent and
$$
D_i:=\left\{T_i\geq \frac{1}{n}; Z_{\theta_i}=-in^{-\gamma}  \right\}.
$$
That is, $D_i$ is the event that the time spent by $Z$ on the interval $I_i$ between times $H_{-(i-1)n^{-\gamma}}$ and $H_{-in^{-\gamma}}$ is greater than $1/n$ and that, after reaching $-(i-1)n^{-\gamma}$, $Z$ exits the interval $I_i$ through $-in^{-\gamma}$.
We have that
$$
\bigcap_{i=1}^n D_i\subset\left\{\int_{-n^{1/(1+\alpha)}}^0 l(\tau_{-n^{1/(1+\alpha)}},u)\rho(du)\geq 1; Z_s \leq n^{-\gamma} \textrm{ for all } s \leq H_{-n^{1/(1+\alpha)}} \right\}.
$$
Observe that the event $D_i$ depends on the realization of $\rho$ only in the interval $I_i$. Also, the family of intervals $I_i$, $i=1,\dots,n$ is disjoint. That, plus the independence of increments of the measure $\rho$ and the strong Markov property of the Brownian motion (applied at the stopping times $\tau_{-in^{-\gamma}}, i=1,\dots,n-1$) implies that the events $D_i,i=1,\dots,n$ are independent. Hence, we have that
\be\label{splitting2}
\prod_{i=1}^n\Pb(D_i)\leq\Pb\left(\int_{-n^{1/(1+\alpha)}}^0 l(\tau_{-n^{1/(1+\alpha)}},u)\rho(du)\geq 1; Z_s \leq n^{-\gamma} \textrm{ for all } s \leq H_{-n^{1/(1+\alpha)}}\right).
\ee

On the other hand, using the scale invariance of the measure $\rho$ and the scale invariance of the Brownian motion $B$, we can show that
$$
\Pb(D_i)=\Pb\left(\int_{-1}^0 l(\tau_{-1},u)\rho(du)\geq1; H_1\geq H_{-1}\right).
$$
That is, the probability of $D_i$ does not depend on $n$ (and is greater than $0$). That, plus display (\ref{splitting2}), proves Lemma \ref{fromminzetatozeta}.
\end{proof}
Now we use Lemmas \ref{hittingtimesfin} and \ref{fromminzetatozeta} to prove Theorem \ref{main2}. We define
$$
G_1:=\left\{ H_{-(n+1)n^{-\gamma}}\leq1 \right\},
$$
$$I^\ast:=[-(2n+1)n^{-\gamma},-(n+1)n^{-\gamma}),$$
$$
T:=\int_{I^\ast}(l(\tau_{-(2n+1)n^{-\gamma}},u)-l(\tau_{-(n+1)n^{-\gamma}},u))\rho(du)\quad\textrm{and}
$$
$$
G_2:=\left\{ T\geq1;Z_s\leq-n^{1/(1+\alpha)}\textrm{ for all } s\in \left[H_{-(n+1)n^{-\gamma}},H_{-(2n+1)n^{-\gamma}} \right]\right\}.
$$

 In the event $G_1$, we have that $Z$ reaches $-(n+1)n^{-\gamma}$ before time $t=1$. On the other hand, in the event $G_2$ we have that, after reaching $-(n+1)n^{-\gamma}$, the process $Z$ remains inside the interval $[-(2n+1)n^{-\gamma},-n^{\alpha/(1+\alpha)}]$ for at least a unit of time. Hence, in the event $G_1\cap G_2$, we have that $Z_1\in [-(2n+1)n^{-\gamma},-n^{\alpha/(1+\alpha)}]$, in particular
$$
 G_1\cap G_2\subset\left\{ Z_1\leq -n^{1/(1+\alpha)}\right\}.
$$
We know, by Lemma \ref{hittingtimesfin}, that $\Pb( H_{-n^{1/(1+\alpha)}}\leq 1 )\geq C_1\exp(-c_1 n)$, for suitable $C_1, c_1>0$. From that, it follows that there exists constants $C,c\geq0$ such  that $\Pb(G_1)\geq C\exp(-cn)$ (to show that, take $M\in\N$ such that $(Mn)^{1/(1+\alpha)}\geq(n+1)n^{-\gamma}$ for all $n\in\N$ and then apply Lemma \ref{hittingtimesfin} at $Mn$). Also, Lemma \ref{fromminzetatozeta} (plus the strong Markov property of the Brownian motion $B$ applied at time $\tau_{-(n+1)n^{-\gamma}}$) implies that $\Pb(G_2)\geq C_2\exp(-c_2)$ for suitable $C_2,c_2>0$. Now, observe that the event $G_1$ depends on the realization of $\rho$ only on the interval $\left[-(n+1)n^{-\gamma},\infty\right)$, whereas $G_2$ depends on the realization of $\rho$ only on the interval $I^\ast$. Those intervals are disjoint, hence, by virtue of the independence of increments of the measure $\rho$ and the strong Markov property of the Brownian motion $B$ applied at the stopping time $\tau_{-(n+1)n^{-\gamma}}$, we deduce that the events $G_1$ and $G_2$ are independent. We have showed that
$$
\Pb\left( |Z_1|\geq n^{1/(1+\alpha)}\right)\geq CC_2\exp(-(c+c_2)n),
$$
and that proves Theorem \ref{main2}.

\end{subsection}
\begin{subsection}{Proof of Theorem \ref{main}}\label{ss:proofofmain1}
The strategy to prove Theorem \ref{main} will be to mimic the arguments leading to Theorem \ref{main2} using the fact that the scaling limit of the one-dimensional BTM is the FIN singular diffusion.
The main tool used in \cite{fin} to prove that the FIN singular diffusion is the scaling limit of the BTM is a coupling between different time scales of the BTM.
We will make use of this coupling for the proof of Theorem \ref{main}, so we proceed to recall it.

To better understand the idea of the coupling, first we have to stress the fact that some discrete-space processes can be expressed as speed-measure changed Brownian motions (see \cite{stone}). In particular, we will express the scalings of the BTM, $(\epsilon X_{\epsilon^{-(1+\alpha)/\alpha}t})_{t\geq0},\epsilon>0$ as speed-measure changed Brownian motions, using a family of random, discrete measures $\rho^\epsilon,\epsilon>0$ as their speed-measures. In fact, we will have that for each $\epsilon>0$, $(B[\rho^\epsilon]_t)_{t\geq 0}$ will be distributed as $(\epsilon X_{\epsilon^{-(1+\alpha)/\alpha}t})_{t\geq0}$, where $B$ is a Brownian motion independent of the $\rho^\epsilon,\epsilon>0$ and $B[\rho^\epsilon]$ is as in \eqref{speedmeasurechangedbrownianmotion}.

On the other hand, the random measures $\rho^\epsilon,\epsilon>0$ will be coupled in such a way that $\rho^\epsilon\stackrel{\epsilon\to0}{\to}\rho$ almost surely, where $\rho$ is the random measure appearing in the definition of the FIN singular diffusion. At this point, we recall that convergence of speed measure changed Brownian motions follows from the convergence of their respective speed measures. More precisely \cite[Theorem 1]{stone}, states that, if $\mu^{\epsilon}$ is a family of measures on $\R$ which converges to some measure $\mu$ as $\epsilon\to 0$, then we have that $(B[\mu^\epsilon]_t)_{t\geq0}$ converges to $(B[\mu]_t)_{t\geq0}$ as $\epsilon\to0$ in the Skorohod topology. Hence, we will have that the convergence of $\rho^\epsilon$ to $\rho$ will imply the convergence of $B[\rho^\epsilon]$ to $B[\rho]$, i.e, we will have that the FIN diffusion is the scaling limit of the BTM.

 Next, we present the construction of such measures.
Let $G:[0,\infty)\to[0,\infty)$ be the function defined by the relation
$$
\Pb(V_1>G(u)):=\Pb(\tau_0>u),
$$
where we recall from the introduction that $(V_x)_{x\in\R}$ is a two sided $\alpha$-stable subordinator and $\tau_0$ is the depth of the trap at $x=0$. The function $G$ is well defined since $V_1$ has a continuous distribution function.
Moreover, $G$ is non-decreasing and right-continuous. Thus, $G$ has a right-continuous generalized inverse $G^{-1}(s):=\inf\{t:G(t)> s\}.$
Now, for all $\epsilon>0$ and $z\in\Z$, we define the random variables $\tau_z^{\epsilon}$ as
$$
\tau^{\epsilon}_z:= G^{-1}(\epsilon^{-1/\alpha}\rho(\epsilon z, \epsilon(z+1)]).
$$
  As $V$ has independent and stationary increments, we have that $(\tau_z^{\epsilon})_{z\in\Z}$ is i.i.d. Moreover, the function $G$ was defined in order to have that $\tau_0^\epsilon$ is distributed according to $\tau_0$. For a proof of that fact we refer to \cite[Proposition 3.1]{fin}.

Finally, we define the coupled family of random measures as
$$
\rho^{\epsilon}:=\sum_{z\in\Z}\epsilon^{1/\alpha}\tau_z^{\epsilon}\delta_{\epsilon z}.
$$
for all $\epsilon>0$.
\begin{obs}\label{rk:measures}
 Define, for each measure $\mu$ over $\R$ and $r\in\R_+$, the following rescaling
\be\label{eq:rescalingofmeasures}
\mu_r(\cdot):=r^{1/\alpha}\mu(r^{-1}\cdot).
\ee
Then $\rho^1_r=\sum_{z\in\Z}r^{1/\alpha} \tau_z^{1}\delta_{r z}$ and, since $(\tau_z^{1})_{z\in\Z}$ is distributed as $(\tau_z^{r})_{z\in\Z}$, we have that
$\rho^1_r$ is distributed as $\rho^r$.
\end{obs}
The next proposition states the two key properties of the coupled measures $\rho^\epsilon, \epsilon>0$.
\begin{prop}\label{lemmafin}(Fontes, Isopi, Newman)
For all $\epsilon>0$ the process $(\epsilon X_{t\epsilon^{-(1+\alpha)/\alpha}})_{t\geq0}$ has the same distribution as $(B[\rho^{\epsilon}]_t)_{t\geq0}$.
Moreover, we have that
\be\label{vagueconvergence}
\rho^{\epsilon}\stackrel{v}{\to}\rho\textrm{   }\Pb\textrm{-a.s.}\textrm{   as }\epsilon\to 0
\ee
where $\stackrel{v}{\to}$ denotes vague convergence of measures.
\end{prop}
\non For the proof of this statement we refer to \cite[Proposition 3.1]{fin}.
Having recalled the coupling between different time scales of the BTM, we turn our attention to the proof of Theorem \ref{main}.

First, note that Proposition \ref{lemmafin} implies, in particular, that $(X_t)_{t\geq0}$ is distributed as $(B[\rho^1]_t)_{t\geq 0}$. Hence, we have that
$$
\Pb(|X_t|\geq x)=\Pb(|B[\rho^1]_t|\geq x)
$$
for all $x\geq0$ and $t \geq 0$.

 For any $b\in\Z$, we define $H^{0}_b:=\inf\{t\geq0: B[\rho^{1}]_t=b\}$. Theorem \ref{main} will be deduced from the two following lemmas, which are the analogous of lemmas \ref{hittingtimesfin} and \ref{fromminzetatozeta} respectively.
\begin{lem}\label{hittingtimesbtm}
There exists positive constants $C,c,\epsilon>0$ such that
$$
\Pb\left(\min\left\{X_s:s\in\left[0,\frac{m^{1/\gamma}}{n^{1/\alpha}}\right]\right\}\leq -m\right)\geq C\exp \left(-cn\right),
$$
for all $m,n\in\N$ such that $n/m\leq\epsilon$.
\end{lem}

\begin{lem}\label{fromminxtox}
 There exists positive constants $C,c,\epsilon>0$ such that
 $$
\Pb\left(\int_{-m}^0l(\tau_{-m},u)\rho^1(du)\geq \frac{m^{1/\gamma}}{n^{1/\alpha}}; Z_s \leq \frac{m}{n} \textrm{ for all } s \leq H^0_{-m} \right)\geq C\exp \left(-c n \right),
$$
for all $m,n\in\N$ such that $n/m\leq\epsilon$.
\end{lem}
First, we will prove Lemma \ref{hittingtimesbtm} by imitating the proof of Lemma \ref{hittingtimesfin}.
Note that, since $X$ is distributed as $B[\rho^1]$, we have that $\{\min\{X_s:s\in[0,\frac{m^{1/\gamma}}{n^{1/\alpha}}]\}\leq -m\}$ is equivalent to $\{H^0_{-m}\leq \frac{m^{1/\gamma}}{n^{1/\alpha}}\}$. Next, define the events $A^0_i, i=1,\dots,n-1$ as
$$
A^0_i:=\left\{\int_{-m\frac{i}{n}}^{-m\frac{i-1}{n}}l(\tau_{-m},u)\rho^1(du)\leq \left(\frac{m}{n}\right)^{1/\gamma}\right\}.
$$
 Also define
$$
A_n^0:=\left\{\int_{-m}^{-m\frac{n-1}{n}}l(\tau_{-m},u)\rho^1(du)\leq \frac{1}{2}\left(\frac{m}{n}\right)^{1/\gamma}\right\}
$$
and
$$
A_0^0:=\left\{\int_{\R_+}l(\tau_{-m},u)\rho^1(du)\leq \frac{1}{2}\left(\frac{m}{n}\right)^{1/\gamma} \right\}.
$$
Clearly, we have that
$$
\bigcap_{i=0}^{n}A^0_i\subset\left\{H^0_{-m}\leq \frac{m^{1/\gamma}}{n^{1/\alpha}}\right\}.
$$
As in the proof of Lemma \ref{hittingtimesfin}, we would like to apply the product rule to the intersection of events above, but those events are not independent. We will deal with the dependence by using the same reasoning as in the proof of Lemma \ref{hittingtimesfin}.
For $i=1,\dots,n$, we define
$$
L_i^0:=\left\{l(\tau_{-m},-m(i-1)/n)\leq m/n\right\}.
$$
The same argument leading to display (\ref{keyestimate}) gives that
\be\label{keyestimatebtm}
\begin{split}
\Pb\left(H^0_{-m}\leq \frac{m^{1/\gamma}}{n^{1/\alpha}} \right)\\
  \geq\Pb\left(A^0_n\cap L^0_n\right)\Pb\left( A^0_0\;\middle\vert\;l(\tau_{-m},0)= m/n \right)\prod_{i=1}^{n-1}&\quad\Pb\left(A^0_i\cap L^0_i\;\middle\vert\;l(\tau_{-m},-im/n)= m/n\right)
 \end{split}
\ee

 Having dealt with the dependence of the events $A^0_i,i=0,\dots,n$, it remains to show that the quantities above can be bounded from below, uniformly on $i,n$ and $m$.
We will adapt the argument in the proof of Lemma \ref{hittingtimesfin}, by using the convergence of $\rho^\epsilon$ to $\rho$.
We start by noting that, for $i=1,\dots,n-1$, $\Pb\left(A^0_i\cap L^0_i\;\middle\vert\;l(\tau_{-m},-im/n)= m/n\right)$ equals
\be\label{nbtm}
\Pb\left(\int_{-m/n}^0 l(\tau_{-m},u)\rho^1\left(du-m\frac{i-1}{n}\right)\leq \left(\frac{m}{n}\right)^{1/\gamma}; l(\tau_{-m},0)\leq m/n \;\middle\vert\; l(\tau_{-m},-m/n)=m/n \right).
\ee
Note that the measure $\rho^1$ is translation invariant in the sense that $\rho^1(\cdot-k)$ is distributed as $\rho^1$ for any $k\in\Z$. Hence, we can replace $\rho^1\left(du-m\frac{i-1}{n}\right)$ by $\rho^{(i)}(du):=\rho^1(du-(m\frac{i-1}{n}-\lfloor m\frac{i-1}{n}\rfloor))$ on display (\ref{nbtm}). Moreover, performing a change of variables inside the integral (and recalling the notation introduced in \eqref{eq:rescalingofmeasures}), we obtain that display (\ref{nbtm}) equals
$$
\Pb\left(\int_{-1}^0 l(\tau_{-m},um/n) \left(\frac{m}{n}\right)^{1/\alpha} \rho^{(i)}_{\frac{n}{m}}(du)\leq \left(\frac{m}{n}\right)^{1/\gamma}; l(\tau_{-m},0)\leq m/n \;\middle\vert\; l(\tau_{-m},-m/n)=m/n \right).
$$
Using the scale invariance of the local time we obtain that the display above is equal to
$$
\Pb\left(\int_{-1}^0 l(\tau_{-m},u) \rho^{(i)}_{\frac{n}{m}}(du)\leq 1; l(\tau_{-m},0)\leq 1 \;\middle\vert\; l(\tau_{-m},-1)=1 \right).
$$
Using independence between the random measure $\rho^1$ and the Brownian motion $B$ we can see that the previous display is at least
\begin{align}\label{eq:independenceoflocaltimeandrho1}
\Pb\left(\sup_{u\in[-1,0]}l(\tau_{-m},u)\leq \delta^{-1}; l(\tau_{-m},0)\leq 1 \;\middle\vert\; l(\tau_{-m},-1)=1 \right)\cdot \Pb\left(\rho^{(i)}_{\frac{n}{m}}(-1,0]\leq \delta\right)
\end{align}
for any $\delta>0$.

We will show that the display above is uniformly bounded away from zero for $\delta<1$ and $m,n\in\N$ with $n/m$ is small enough.
It follows from standard considerations about Bessel processes that the first factor in \eqref{eq:independenceoflocaltimeandrho1} is positive for $\delta<1$. Then, it suffices to bound from below $\Pb\left(\rho^{(i)}_{\frac{n}{m}}(-1,0]\leq \delta\right)$.
At this point we will make use of the convergence of $\rho^{\epsilon}\to\rho$ as $\epsilon\to0$. From Observation \ref{rk:measures} it follows that $\rho^{(i)}_{\frac{n}{m}}$ is distributed as $\rho^{n/m}(\cdot-\frac{n}{m}(m\frac{i-1}{n}-\lfloor m\frac{i-1}{n}\rfloor))$. On the other hand, using display (\ref{vagueconvergence}) in Lemma \ref{lemmafin} and the fact that $m\frac{i-1}{n}-\left\lfloor m\frac{i-1}{n}\right\rfloor\leq1$ we obtain that $\rho^{(i)}_{\frac{n}{m}}=\rho^{n/m}\left(du-\frac{n}{m}\left(m\frac{i-1}{n}-\left\lfloor m\frac{i-1}{n}\right\rfloor\right)\right)$ converges to $\rho$ as $\frac{n}{m}$ goes to $0$. Moreover, we see that $\Pb(\rho^{(i)}_{\frac{n}{m}}(-1,0]\leq \delta)$ converges to $\Pb(\rho(-1,0]\leq\delta)$ uniformly on $i$ as $\frac{n}{m}$ goes to $0$, for each $\delta >0$.
Thus, there exists $\epsilon$ small enough such that, for $n,m\in\N$ such that $n/m\leq\epsilon$ we have that the quantity on display (\ref{nbtm}) is bounded away from zero uniformly on $n,m$ and $i$. Following analogous arguments we can find uniform lower bounds for $\Pb\left(A^0_n\cap L^0_n\right)$ and $\Pb\left( A^0_0\;\middle\vert\;l(\tau_{-x},0)= m/n \right)$. That, plus display (\ref{keyestimatebtm}), proves Lemma \ref{hittingtimesbtm}.

In a similar way, we can apply the vague convergence of $\rho^\epsilon$ to $\rho$ to obtain Lemma \ref{fromminxtox} by imitating the argument leading to Lemma \ref{fromminzetatozeta}.

Now we are ready to prove Theorem \ref{main}. As both, $\Pb(|X_t|\geq x)$ and $C\exp(-c(\frac{x}{t^\gamma})^{1+\alpha})$ are decreasing on $x$, to prove Theorem \ref{main}, it will suffice to show that there exist constants $C,c,\epsilon_1>0$ such that
$$
\Pb(|X_t|\geq m)\geq C\exp \left(-c \left(\frac{m}{t^\gamma}\right)^{1+\alpha}\right)
$$
for all $m\in\N$ and $t$ such that $\frac{m}{\epsilon}\leq t$. Moreover, we can restrict ourselves to $t\leq m^{1/\gamma}$, because for $t\geq m^{1/\gamma}$, we have that $\Pb(|X_t|\geq m)\geq\Pb(|X_t|\geq t^\gamma)$, and the convergence of $(\epsilon^{-1}X_{\epsilon^{1/\gamma}t})_{t\in\R_+}$ to the FIN singular diffusion (proved in \cite[Theorem 4.1]{fin}) implies that $\Pb(|X_t|\geq t^\gamma)$ converges to $\Pb(|Z_1|\geq 1)$ as $t\to\infty$.

Let $n\geq2$ be a natural number such that $t\in\left[\frac{m^{1/\gamma}}{n^{1/\alpha}},\frac{m^{1/\gamma}}{(n-1)^{1/\alpha}}\right]$ (which exists because $t\leq m^{1/\gamma}$).
Let $K\in\N$ be a fixed constant that satisfies $K^{\alpha}>2$. Let us define
$$
G_1^0:=\left\{\min\left\{X_s:s\in\left[0,\frac{m^{1/\gamma}}{n^{1/\alpha}}\right]\right\}\leq-\frac{n+1}{n}m\right\}
$$
$$
T^0:=\int_{-\frac{(K+1)n+1}{n}m}^{-\frac{n+1}{n}m}(l(\tau_{-\frac{(K+1)n+1}{n}m},u)-l(\tau_{-\frac{n+1}{n}m},u))\rho^1(du),
$$
and
$$
G_2^0:=\left\{ T^0\geq\frac{(Km)^{1/\gamma}}{(Kn)^{1/\alpha}};X_s\leq-m\textrm{ for all } s\in \left[H^0_{-\frac{n+1}{n}m},H^0_{-\frac{(K+1)n+1}{n}m} \right]\right\}.
$$
In the event $G_1^0$, $X$ reaches $-\frac{n+1}{n}m$ before time $\frac{m^{1/\gamma}}{n^{1/\alpha}}$ and, since we have chosen $n$ such that $t\geq\frac{m^{1/\gamma}}{n^{1/\alpha}}$,
this implies that $X$ reaches $-\frac{n+1}{n}m$ before time $t$. On the other hand, in the event $G_2^0$ we have that, after reaching $-\frac{n+1}{n}m$, $X$ remains inside the interval $[-\frac{(K+1)n+1}{n}m,-m]$ for a time which is bigger than $\frac{(Km)^{1/\gamma}}{(Kn)^{1/\alpha}}$. But $t\leq\frac{(Km)^{1/\gamma}}{(Kn)^{1/\alpha}}$, because $t\leq\frac{m^{1/\gamma}}{(n-1)^{1/\alpha}}$ and $K^{\alpha}>2\geq \frac{n}{n-1}$. Hence, we have that
$$
 G^0_1\cap G^0_2\subset \left\{ X_t\leq -m\right\}.
$$
As in the proof of Theorem \ref{main2}, we have that $G_1^0$ and $G_2^0$ are independent. Thus we can apply Lemma \ref{hittingtimesbtm} to find that $\Pb(G_1)\geq C_1\exp(-c_1 n)$. We can also apply Lemma \ref{fromminxtox} (at $Km$ and $Kn$) to find that $\Pb(G_2)\geq C_2\Pb(c_2 Kn)$. Thus we have that
$$
\Pb(X_t\leq -m)\geq C_1C_2\exp(-(c_1+Kc_2)n).
$$
as we have chosen $n$ such that $t\leq \frac{m^{1/\gamma}}{(n-1)^{1/\alpha}}$, we have that $n\leq(\frac{m}{t^\gamma})^{1+\alpha}+1$ and we obtain
$$
\Pb(X_t\leq -m)\geq C_1C_2\exp(-(c_1+Kc_2))\exp\left(-(c_1+Kc_2)\left(\frac{m}{t^\gamma}\right)^{1+\alpha}\right).
$$
The last display provides the lower desired lower bound with $C=C_1C_2\exp(-(c_1+Kc_2))$
and $c=c_1+Kc_2$.
\end{subsection}
\end{section}

\bibliographystyle{amsplain}
\bibliography{manolo}
\end{document}